\documentclass{birkmult}  % amslatex

\usepackage{amssymb}
\usepackage{amscd}
\usepackage{amsmath}

\newtheorem{thm}{Theorem}[subsection]

\newtheorem{lemma}[thm]{Lemma}

\newtheorem{definition}[thm]{Definition}

\theoremstyle{remark}
\newtheorem{remark}[thm]{Remark}
\newtheorem{example}[thm]{Example}
\theoremstyle{definition}

\numberwithin{equation}{section}

\newcommand{\cA}{\mathcal{A}}
\newcommand{\cB}{\mathcal{B}}
\newcommand{\cC}{\mathcal{C}}
\newcommand{\cD}{\mathcal{D}}
\newcommand{\cE}{\mathcal{E}}
\newcommand{\cF}{\mathcal{F}}
\newcommand{\cM}{\mathcal{M}}
\newcommand{\cO}{\mathcal{O}}

\newcommand{\bfrho}{\boldsymbol{\rho}}

\DeclareMathOperator{\Perf}{Perf}
\DeclareMathOperator{\CC}{CC}
\DeclareMathOperator{\Matr}{Matr}

\newcommand{\tw}{\mathrm{tw}}

\DeclareMathOperator{\Ad}{Ad} \DeclareMathOperator{\Hom}{Hom}
\DeclareMathOperator{\Fun}{Fun}
\newcommand{\limdir}{\mathop{\varinjlim}\limits}
\newcommand{\liminv}{\mathop{\varprojlim}\limits}

\begin{document}

\title{Chern character for twisted complexes}

\author[P.Bressler]{Paul Bressler}
\address{School of Mathematics, Institute for Advanced Study}
\email{bressler@math.ias.edu}

\author[A.Gorokhovsky]{Alexander Gorokhovsky}
\address{Department of Mathematics, University of Colorado}
\email{Alexander.Gorokhovsky@colorado.edu}

\author[R.Nest]{Ryszard Nest}
\address{Department of Mathematics, University of Copenhagen}
\email{rnest@math.ku.dk}

\author[B.Tsygan]{Boris Tsygan}
\address{Department of Mathematics, Northwestern University}
\email{tsygan@math.northwestern.edu}

\maketitle

\begin{center}
In memory of Sasha Reznikov
\end{center}

\bigskip
\section{Introduction}
The Chern character from the algebraic K theory to the cyclic
homology of associative algebras was defined by Connes and Karoubi
\cite{C}, \cite{K}, \cite{L}. Goodwillie and Jones \cite{Go}, \cite{J} defined the
negative cyclic homology and the Chern character with values there.
In this paper we generalize this Chern character to the K theory of
twisted modules over twisted sheaves of algebras.

More precisely, we outline the construction of the Chern character
of a perfect complex of twisted sheaves of modules over an algebroid
stack $\mathcal{A}$ on a space $M$. This includes the case of a
perfect complex of sheaves of modules over a sheaf of algebras
$\mathcal{A}$. In the latter case, the recipient of the Chern
character is the hypercohomology of $M$ with coefficients in the
sheafification of the presheaf of negative cyclic complexes. The
construction of the Chern character for this case was given in
\cite{BNT1} and \cite{K}. In the twisted case, it is not a priori
clear what the recipient should be. One can construct \cite{K2}, \cite{MC}  the Chern character with values in the negative cyclic homology
of the category of perfect complexes (localized by the subcategory of acyclic complexes); the question is, how to
compute this cyclic homology, or perhaps how to map it into
something simpler.

Ideally, the recipient of the Chern character would be the
hypercohomology of $M$ with coefficients in the negative cyclic
complex of a sheaf of associative algebras. We show that
this is almost the case. We construct associative algebras that form
a presheaf not exactly on $M$ but rather on a first barycentric
subdivision of the nerve of a cover of $M$. These algebras are
twisted matrix algebras. We used them in \cite{BGNT} and
\cite{BGNT1} to classify deformations of algebroid stacks.

We construct the Chern characters
\begin{eqnarray}
K_\bullet(\Perf(\cA)) & \xrightarrow{} &
\check{\mathbb{H}}^{-\bullet}(M, \CC^-_{\bullet}(\Matr_\tw(\cA)))
\label{eq:ch0} \\
K_\bullet(\Perf_Z(\cA)) & \xrightarrow{} &
\check{{\mathbb H}}_Z^{-\bullet}(M,
{\CC}^-_{\bullet}(\Matr_\tw(\cA)))
\label{eq:ch0 Z}
\end{eqnarray}
where $K_\bullet(\Perf(\cA))$ is the $K$ theory of perfect complexes
of twisted $\cA$-modules, $K_\bullet(\Perf_Z(\cA))$
is the $K$ theory of perfect complexes of twisted $\cA$-modules
acyclic outside a closed subset $Z$, and the right hand sides are
the hypercohomology of $M$ with coefficients in the negative cyclic
complex of twisted matrices, cf. Definition \ref{dfn:complex of cc
of matr} .

Our construction of the Chern character is more along the lines of
\cite{K} than of \cite{BNT1}. It is modified for the twisted case
and for the use of twisted matrices. Another difference is a method
that we use to pass from perfect to very strictly perfect complexes. This
method involves a general construction of operations on cyclic
complexes of algebras and categories. This general construction, in
partial cases, was used before in \cite{NT}, \cite{NT1} as a version
of noncommutative calculus. We recently realized that it can be
obtained in large generality by applying the functor
${\CC}^-_{\bullet}$ to the categories of $A_\infty$
functors from \cite{BLM}, \cite{K1}, \cite{Ko}, \cite{Lu}, and
\cite{Ta}.

The fact that these methods are applicable is due to the observation
that a perfect complex, via the formalism of twisting cochains of
O'Brian, Toledo, and Tong, can be naturally interpreted as an
$A_\infty$ functor from the category associated to a cover to the
category of strictly perfect complexes. The fourth  author is
grateful to David Nadler for explaining this to him.

In the case when the stack in question is a gerbe, the recipient of the Chern character maps to the De Rham cohomology twisted by the three-cohomology class determined by this gerbe (the Dixmier-Douady class). A Chern character with values in the twisted cohomology was constructed in \cite{MaS}, \cite{BCMMS}, \cite{AS} and generalized in \cite{MaS1} and \cite{TX}.  The $K$-theory which is the source of this Chern character is rather different from the one studied here. It is called the twisted $K$-theory and is a generalization of the topological $K$-theory. Our Chern character has as its source the algebraic $K$-theory which probably maps to the topological one.

\vskip 5mm
\noindent\textbf{Acknowledgements:} P.B. was supported by the Ellentuck Fund, A.G. was supported by supported by NSF grant DMS-0400342, B.T. was supported by NSF grant
DMS-0605030.

\section{Gerbes and stacks}\label{Gerbes and stacks}

\subsection{}
Let $M$ be a topological space. In this paper,
by a stack on $M$ we will mean an equivalence class of the following
data:
\begin{enumerate}
\item an open cover $M=\cup U_i$;

\item a sheaf of rings $\cA _i$ on every $U_i$;

\item an isomorphism of sheaves of rings $G_{ij}: \cA _j|(U_i \cap
U_j) \cong \cA _i |(U_i \cap U_j)$ for every $i,\,j$;

\item an invertible element $c_{ijk} \in \cA _i (U_i \cap U_j \cap
U_k)$ for every $i,\,j,\,k$ satisfying
\begin{equation} \label {eq:2-cocycle 1}
G_{ij}G_{jk}=\Ad(c_{ijk})G_{ik}
\end{equation}
such that, for every $i,\,j,\,k,\,l$,
\begin{equation} \label {eq:2-cocycle 2}
c_{ijk}c_{ikl}=G_{ij}(c_{jkl})c_{ijl}
\end{equation}
\end{enumerate}

To define equivalence, first recall the definition of a refinement.
An open cover $\mathfrak{V}=\{V_j\}_{j\in J}$ is a refinement of an
open cover $\mathfrak{U}=\{U_i\}_{i\in I}$ if a map $f: J\to I$ is
given, such that $V_j \subset U_{f(j)}.$ Open covers form a
category: to say that there is a morphism from $\mathfrak{U}$ to
$\mathfrak{V}$ is the same as to say that $\mathfrak{V}$ is a
refinement of $\mathfrak{U}$. Composition corresponds to composition
of maps $f$.

Our equivalence relation is by definition the weakest for which the
two data $(\{U_i\}, \cA_i, G_{ij}, c_{ijk})$ and
\[
(\{V_p\}, \cA_{f(p)}|V_p, G _{f(p)f(q)}, c_{f(p)f(q)f(r)})
\]
are equivalent whenever $\{V_p\}$ is a refinement of $\{U_i\}$ (the
corresponding map $\{p\}\to \{i\}$ being denoted by $f$).

If two data $(\{U'_i\},\; \cA' _i, \; G'_{ij},\; c'_{ijk})$ and
$(\{U''_i\},\; \cA'' _i, \; G''_{ij},\; c''_{ijk})$ are given on
$M$, define an isomorphism between them as follows. First, choose an
open cover  $M=\cup U_i$ refining both $\{U'_i\}$ and $\{U''_i\}$.
Pass from our data to new, equivalent data corresponding to this
open cover. An isomorphism is  an equivalence class of a collection
of isomorphisms $H_i: \cA' _i \cong \cA''_i$ on $U_i$ and invertible
elements $b_{ij}$ of $\cA' _i (U_i \cap U_j)$ such that
\begin{equation} \label{eq:equivalence of stacks 1}
G''_{ij}=H_i \Ad(b_{ij})G'_{ij}H_j ^{-1}
\end{equation}
and
\begin{equation} \label{eq:equivalence of stacks 2}
H_i^{-1}(c''_{ijk})=b_{ij}G'_{ij}(b_{jk})c'_{ijk}b_{ik}^{-1}
\end{equation}

If $\{V_p\}$ is a refinement of $\{U_i\}$, we pass from $(\{U_i\}, \cA_i,
G_{ij}, c_{ijk})$ to the equivalent data $(\{V_p\}, A _{f(p)f(q)},
c_{f(p)f(q)f(r)})$ as above. We define the equivalence relation to
be the weakest for which, for all refinements, the data $(H_i,
b_{ij})$ and $(H_{f(p)}, b_{f(p)f(q)})$ are equivalent.

Define composition of isomorphisms as follows. Choose a common
refinement $\{U_i\}$ of the covers $\{U'_i\},$ $\{U''_i\},$ and
$\{U'''_i\}.$ Using the equivalence relation, identify all the stack
data and all the isomorphism data with the data corresponding to the
cover $\{U_i\}$.  Define $H_i = H'_i\circ H''_i$ and
$b_{ij}={H''_i}^{-1}(b'_{ij})b''_{ij}$.  It is easy to see that this
composition is associative and is well defined for equivalence
classes.

Now consider two isomorphisms $(H'_i, b'_{ij})$ and $(H''_i,
b''_{ij})$ between the stacks $(\{U'_i\},\; \cA' _i, \; G'_{ij},\;
c'_{ijk})$ and $(\{U''_i\},\; \cA'' _i, \; G''_{ij},\; c''_{ijk})$.
We can pass to a common refinement, replace our data by equivalent
data, and assume that $\{U'_i\}=\{U''_i\}=\{U_i\}.$ A two-morphism
between the above isomorphisms is an equivalence class of a
collection of invertible elements $a_i$ of $\cA'_i(U_i)$ such that
$H_i''=H_i'\circ {\operatorname{Ad}}(a_i)$ and
$b''_{ij}=a_i^{-1}b'_{ij}G'_{ij}(a_j).$ The equivalence relation is
the weakest for which, whenever $\{V_p\}$ is a refinement of
$\{U_i\}$, $\{a_i\}$ is equivalent to $\{a_{f(p)}\}:(H'_{f(p)},
b'_{f(p)f(q)})\to (H''_{f(p)}, b''_{f(p)f(q)}) $. The composition
between $\{a'_i\}$ and $\{a''_i\}$ is defined by $a_i=a'_i a''_i.$
This operation is well-defined  at the level of equivalence classes.

With the operations thus defined, stacks form a two-groupoid.

\emph{A gerbe} on a manifold $M$ is a stack for which $\cA _i$ = $\cO _{U_i}$ and
$G_{ij}=1$. Gerbes are classified up to isomorphism by cohomology
classes in $H^2 (M, \cO^*_M).$

For a stack $\cA$ define \emph{a twisted $\cA$-module} over an open
subset $U$ as an equivalence class of a collection of sheaves of
$\cA_i$ -modules $\cM _i$ on $U\cap U_i$, together with isomorphisms
$g_{ij}: \cM _j \to \cM _i$ on $U\cap U_i \cap U_j$ such that
$g_{ik}=g_{ij}G_{ij}(g_{jk})c_{ijk}$ on $U\cap U_i\cap U_j \cap
U_k.$ The equivalence relation is the weakest for which, if
$\{V_p\}$ is a refinement of $\{U_i\}$, the data $(\cM_{f(p)},
g_{f(p)f(q)})$ and $(\cM_i, g_{ij})$ are equivalent.

We leave it to the reader to define morphisms of twisted modules.
A twisted module is said to be \emph{free} if the $\cA_i$-module $\cM_i$ is.

\subsection{Twisting cochains}\label{Twisting cochains}
Here we recall the formalism from \cite{TT}, \cite{OTT}, \cite{OB}, generalized to
the case of stacks. For a stack on $M=\cup U_i$ as above, by $\cF$ we will denote a collection $\{\cF_i\}$ where $\cF_i$ is a graded sheaf which is a direct summand of a free graded
$\cA _i$-module of finite rank on $U_i.$ A $p$-cochain with values in $\cF$ is a collection $a _{i_0 \ldots i_p} \in \cF _{i_0}
(U_{i_0}\cap \ldots \cap U_{i_p}); $ for two collections $\cF$ and $\cF'$ as above, a $p$-cochain with values in
${\operatorname {Hom}}(\cF,\cF')$ is a collection $a _{i_0
\ldots i_p} \in {\operatorname {Hom}}_{\cA _{i_0}}(\cF _{i_p},
\cF '_{i_0}) (U_{i_0}\cap \ldots \cap U_{i_p})$ (the sheaf $\cA
_{i_0}$ acts on $\cF _{i_p}$ via $G_{i_0i_p}$). Define the cup
product by
\begin{equation} \label{eq:cup}
(a\smile b)_{i_0 \ldots i_{p+q}}=(-1)^{|a_{i_0\ldots i_p}|q}a_{i_0\ldots i_p}G_{i_pi_{p+q}}(b_{i_{p+1}\ldots i_{p+q}})c_{i_0i_pi_{p+q}}
\end{equation}
and the differential by
\begin{equation} \label{eq:differential}
(\check{\partial }a)_{i_0 \ldots i_{p+1}}=\sum _{k=1}^p (-1)^{k}a_{i_0\ldots \widehat{i_k} \ldots i_{p+1}}
\end{equation}
Under these operations, ${\operatorname {Hom}}(\cF,\cF)$-valued cochains form a DG algebra and $\cF$-valued cochains a DG module. 

If ${\mathfrak{V}}$ is a refinement of ${\mathfrak{U}}$ then cochains with respect to ${\mathfrak{U}}$ map to cochains with respect to ${\mathfrak{V}}$. For us, the space of cochains will be always understood as the direct limit over all the covers.

\emph{A twisting cochain} is a $\Hom(\cF,\cF)$-valued cochain $\rho$ of total degree one such that
\begin{equation}\label{eq:twisting cochain}
\check{\partial }\rho +\frac{1}{2}\rho \smile \rho =0
\end{equation}
A morphism between twisting cochains $\rho$ and $\rho '$ is a cochain $f$ of total degree zero such that $\check{\partial }f +\rho'\smile f-f\smile \rho = 0$. A homotopy between two such morphisms $f$ and $f'$ is a cochain $\theta$ of total degree $-1$ such that $f-f'={\check{\partial}}\theta+\rho' \smile \theta + \theta \smile \rho .$ More generally, twisting cochains form a DG category. The complex $\Hom(\rho, \rho ')$ is the complex of $\Hom(\cF,\cF')$-valued cochains with the differential
$$
f\mapsto {\check{\partial}}f+\rho'\smile f-(-1)^{|f|}f\smile \rho \ .
$$

There is another, equivalent definition of twisting cochains. Start
with a collection $\cF=\{\cF_i\}$ of direct summands of free graded twisted
modules of finite rank on $U_i$ (a twisted module on $U_i$ is said to be free if the corresponding $\cA_i$-module is). Define ${\operatorname{Hom}}(\cF,\cF')$-valued cochains as collections
of morphisms of graded twisted modules $a_{i_0\ldots
i_p}:\cF_{i_p}\to \cF'_{i_0}$ on $ U_{i_0}\cap \ldots \cap U_{i_p}$.
The cup product is defined by
\begin{equation} \label{eq:cup 1}
(a\smile b)_{i_0 \ldots i_{p+q}}=(-1)^{|a_{i_0\ldots i_p}|q}a_{i_0\ldots i_p}b_{i_{p+1}\ldots i_{p+q}}
\end{equation}
and the differential by \eqref{eq:differential} . A twisting cochain is a cochain $\rho$ of total degree $1$ satisfying \eqref{eq:twisting cochain}.

If one drops the requirement that the complexes $\cF$ be direct
summands of graded free modules of finite rank, we get objects that
we will call \emph{weak twisting cochains}. A morphism of (weak)
twisting cochains is \emph{a quasi-isomorphism} if $f_i$ is for
every $i$. Every complex $\cM$ of twisted modules can be viewed as a
weak twisting cochain, with $\cF_i=\cM$ for all $i$,
$\rho_{ij}=\operatorname{id}$ for all $i,\,j$, $\rho_i$ is the
differential in $\cM,$ and $\rho_{i_0\ldots i_p}=0$ for $p> 2$. We
denote this weak twisting cochain by $\rho_0({\cM}).$ By $\bfrho_0$
we denote the DG functor $\cM\mapsto \rho_0(\cM)$ from the DG
category of perfect complexes to the DG category of weak twisting
cochains.

If $\{V_s\}$ is a refinement of $\{U_i\}$, we declare twisting cochains $(\cF _i, \rho _{i_0\ldots i_p})$ and $(\cF_{f(s)}|V_s, \rho _{f(s_0)\ldots f(s_p)})$ equivalent. Similarly for morphisms.

A complex of twisted modules is called \emph{perfect} (resp. {\em strictly perfect)} if it is locally isomorphic in the derived category (resp. isomorphic) to a direct summand of a bounded complex of finitely generated free modules. A parallel definition can of course be given for complexes of modules over associative algebras.
\begin{lemma}\label{lemma:perfect complex} Let $M$ be paracompact.

1. For a perfect complex $\cM$ there exists a twisting cochain $\rho$ together with a quasi-isomorphism of weak twisting cochains
$\rho \xrightarrow{\phi} \rho_0({\cM})$.

2. Let $f: \cM_1\to \cM_2$ be a morphism of perfect complexes. Let $\rho _i, \phi _i$ be twisting cochains corresponding to $\cM_i$, $i=1,2.$ Then there is a morphism of twisting cochains $\varphi(f)$ such that $\phi_2 \varphi (f) $ is homotopic to $f\phi_1.$

3. More generally, each choice $\cM\mapsto \rho(\cM)$ extends to an $A_\infty$ functor $\bfrho$ from the DG category of perfect complexes to the DG category of twisting cochains, together with an $A_\infty$ quasi-isomorphism $\bfrho\to \bfrho_0$. (We recall the definition of $A_\infty$ functors in \ref{Hochschild}, and that of $A_\infty$ morphisms of $A_\infty$ functors in \ref{ss:Categories of functors}).
\end{lemma}

{\em Sketch of the proof:} 
We will use the following facts about complexes of modules over associative algebras.

1) If a complex $\cF$
is strictly perfect, for a quasi-isomorphism $\psi:\cM\to \cF$ there
is a quasi-isomorphism $\phi:\cF\to\cM$ such that $\psi\circ\phi$ is
homotopic to the identity. 

2) If $f: \cM_1\to \cM_2$ is a morphism of
perfect complexes and $\phi _i: \cF_i\to \cM_i$ , $i=1,2,$ are
quasi-isomorphisms with $\cF_i$ strictly perfect, then there is a
morphism $\varphi(f):\cF_1\to\cF_2$ such that $\phi_2 \varphi  $ is
homotopic to $f\phi_1.$ 

3) If $\cF$ is strictly perfect and
$\phi: \cF\to\cM$ is a morphism which is zero on cohomology, then
$\phi$ is homotopic to zero. 

Let $\cM$ be a perfect complex of twisted
modules. Recall that, by our definition, locally, there is a chain of quasi-isomorphisms connecting it to a strictly perfect complex $\cF$. Let us start by observing that one can replace that by a quasi-isomorphism from $\cF$ to $\cM$. In other words, locally, there is a strictly perfect complex $\cF$
and a quasi-isomorphism $\phi: \cF\to \cM.$ Indeed, this is true at the level of germs at every point, by virtue of 1) above. For any point, the images of generators of $\cF$ under morphisms $\phi$, resp. under homotopies $s,$ are germs of sections of $\cM,$ resp. of $\cF,$ which are defined on some common neighborhood. Therefore quasi-isomorphisms and homotopies are themselves defined on these neighborhoods.

We get a cover $\{U_i\}$, strictly perfect complexes $\cF_i$ with differentials $\rho_i,$ and quasi-isomorphisms $\phi_i:\cF_i\to \cM$ on $U_i.$ Now observe that, at any point of $U_{ij},$ the morphisms $\rho_{ij}$ can be constructed at the level of germs because of 2). As above, we conclude that each of them can be constructed on some neighborhood of this point. Replace the cover $\{U_i\}$ by a locally finite refinement$\{U'_i\}$. Then, for every point $x,$ find a neighborhood $V_x$ on which all $\rho_{ij}$ can be constructed. Cover $M$ by such neighborhoods. Then pass to a new cover which is a common refinement of $\{U'_i\}$ and $\{V_x\}$. For this cover, the component $\rho_{ij}$ can be defined.

Acting as above, using 2) and 3), one can
construct all the components of the twisting cochain $\rho(\cM)$, of
the $A_\infty$ functor $\bfrho,$ and of the $A_\infty$ morphism of
$A_\infty$ functors $\bfrho\to\bfrho_0.$
\begin{remark}\label{rmk:simplification} One can assume that all components of a twisting cochain $\rho$ lie in the space of cochains with respect to one and the same cover if the following convention is adopted:  all our perfect complexes are locally quasi-isomorphic to strictly perfect complexes as {\em complexes of presheaves}. In other words, there is an open cover $\{U_i\}$ together with a strictly perfect complex $\cF_i$  and a morphism $\phi_i:\cF_i\to \cM$ on any $U_i,$ such that $\phi_i$ is a quasi-isomorphism at the level of sections on any open subset of $U_i.$
\end{remark}
%Let $C^{\bullet}(M, {\operatorname{End}}(\cM))$, resp. $C^{\bullet}(M,\cM),$ be the complex of cochains with values in ${\Hom}(\cF, \cF),$ resp. in $\cF$, for a twisting cochain from the above lemma. This complexes are defined canonically up to quasi-isomorphism.
\subsection{Twisted matrix algebras}
For any $p$-simplex $\sigma$ of the nerve of an open cover $M=\cup U_i$ which corresponds to $U_{i_0}\cap \ldots \cap  U_{i_p}$, put $I_{\sigma}=\{i_0, \ldots, i_p\}$ and $U_{\sigma}=\cap_{i\in {I_\sigma}}U_i.$ Define the algebra ${\rm {Matr}}^{\sigma}_{tw}(\cA)$ whose elements are finite matrices
$$\sum_{i,j\in I_{\sigma}} a_{ij}E_{ij}$$
such that $a_{ij} \in (\cA _i(U_{\sigma})). $
The product is defined by
$$a_{ij}E_{ij}\cdot  a_{lk}E_{lk} = \delta_{jl} a_{ij}G_{ij}(a_{jk})c_{ijk}E_{ik}$$
For $\sigma\subset \tau,$ the inclusion
$$i_{\sigma\tau}:{\rm {Matr}}^{\sigma}_{tw}(\cA)\to {\rm {Matr}}^{\tau}_{tw}(\cA),$$
${\sum} a_{ij}E_{ij}\mapsto{\sum} (a_{ij}|U_{\tau})E_{ij},$ is a morphism of algebras (not of algebras with unit).Clearly, $i_{\tau\rho}i_{\sigma\tau}=i_{\sigma\rho}.$ If ${\mathfrak V}$ is a refinement of ${\mathfrak U}$ then there is a map
$${\rm {Matr}}^{\sigma}_{tw}(\cA)\to {\rm {Matr}}^{f(\sigma)}_{tw}(\cA)$$
which sends ${\sum} a_{ij}E_{ij}$ to ${\sum} (a_{f(i)f(j)}|V_{f(\sigma)})E_{f(i)f(j)}$.

\begin{remark}\label{rmk:Morita}
For a nondecreasing map $f: I_\sigma \to I_\tau$ which is not necessarily an inclusion, we have the bimodule $M_f$ consisting of twisted $|I_\sigma| \times |I_\tau|$ matrices. Tensoring by this bimodule defines the functor
$$f_*: {\rm {Matr}}^{\sigma}_{tw}(\cA)-{\rm{mod}}\to {\rm {Matr}}^{\tau}_{tw}(\cA)-{\rm{mod}}$$
such that $(fg)_*=f_*g_*.$
\end{remark}

\section{The Chern character}\label{ss:The Chern character for rings}
\subsection{Hochschild and cyclic complexes}\label{Hochschild} We start by recalling some facts and constructions from noncommutative geometry. Let $A$ be an associative unital algebra over a unital algebra $k$. Set
%$$\overline{A} = A / k\cdot 1$$
%and
$$C_p (A,A) = C_p(A) = {A} ^{\otimes (p+1)}.$$
We denote by $b:C_p (A)\to C_{p-1} (A)$ and $B:C_p (A)\to C_{p+1} (A)$ the standard differentials from the Hochschild and cyclic homology theory (cf. \cite{C}, ,\cite{L}, \cite{T}). The Hochschild chain complex is by definition $(C_{\bullet} (A),b);$ define
$$CC^-_{\bullet}(A)=(C_{\bullet} (A)[[u]],b+uB);$$
$$CC^{\operatorname{per}}_{\bullet}(A)=(C_{\bullet} (A)[[u,u^{-1}],b+uB);$$
$$CC_{\bullet}(A)=(C_{\bullet} (A)[[u,u^{-1}]/uC_{\bullet} (A)[[u]],b+uB).$$
These are, respectively, {\em the negative cyclic}, {\em the periodic cyclic}, and {\em the cyclic} complexes of $A$ over $k$.

We can replace $A$ by a small DG category or, more generally, by a small $A_\infty$ category. Recall that a small $A_\infty$ category consists of a set $\mathrm{Ob}(\cC)$ of objects and a graded $k$-module of ${\cC}(i,j)$ of morphisms for any two objects $i$ and $j$, together with compositions
$$m_n: {\cC}(i_{n},i_{n-1})\otimes \ldots \otimes {\cC}(i_{1},i_0)\to {\cC}(i_{n},i_0)$$
of degree $2-n,$ $n\geq 1$, satisfying standard quadratic relations to which we refer as the $A_\infty$ relations. In particular, $m_1$ is a differential on ${\cC}(i,j)$. An $A_\infty$ functor $F$ between two small $A_\infty$ categories $\cC$ and $\cD$ consists of a map $F: {\rm{Ob}}(\cC)\to {\rm{Ob}}(\cD)$ and $k$-linear maps
$$F_n: {\cC}(i_{n},i_{n-1})\otimes \ldots \otimes{\cC}(i_{1},i_0)\to {\cD}(Fi_{n},Fi_0)$$
of degree $1-n,$ $n\geq 1,$
satisfying another standard relation. We refer the reader to \cite{K1} for formulas and their explanations.

For a small $A_\infty$ category $\cC$ one defines the Hochschild complex $C_{\bullet}(\cC)$ as follows:
$$C_{\bullet}(\cC)=\bigoplus _{i_0,\ldots,i_n\in {\rm{Ob}}(\cC)}{\cC}(i_{1},i_{0})\otimes {\cC}(i_{2},i_{1})\otimes\ldots \otimes{\cC}(i_{n},i_{n-1})\otimes {\cC}(i_{0},i_n)$$
(the total cohomological degree being the degree induced from the grading of $\cC(i,j)$ minus $n$). The differential $b$ is defined by
$$b(f_0\otimes\ldots f_n)=\sum_{j,k}\pm m_k(f_{n-j+1}, \ldots, f_0, \ldots, f_{k-1-j})\otimes f_{k-j}\otimes \ldots \otimes f_{n-j}$$
$$+\sum_{j,k}\pm f_0\otimes \ldots \otimes f_j \otimes m_k(f_{j+1},\ldots, f_{j+k})\otimes \ldots \otimes f_n$$
The cyclic differential $B$ is defined by the standard formula with appropriate signs; cf. \cite{G}.
\subsection{Categories of $A_\infty$ functors} \label {ss:Categories of functors}

For two $DG$ categories $\cC$ and $\cD$ one can define the DG category $\Fun_\infty(\cC,\cD)$. Objects of $\Fun_\infty(\cC,\cD)$ are $A_\infty$ functors $\cC\to \cD$. The complex $\Fun_\infty(\cC,\cD)(F,G)$ of morphisms from $F$ to $G$ is the Hochschild cochain complex of $\cC$ with coefficients in $\cD$ viewed as an $A_\infty$ bimodule over $\cC$ via the $A_\infty$ functors $F$ and $G$, namely
$$\prod _{i_0,\ldots,i_n\in {\rm{Ob}}(\cC)}{\Hom}(\cC(i_0,i_1)\otimes \ldots \otimes \cC(i_{n-1},i_n), \cD(Fi_0,Gi_n))$$
The DG category structure on $\Fun_\infty(\cC, \cD)$ comes from the cup product. More generally, for two $A_\infty$ categories $\cC$ and $\cD$, $\Fun_\infty(\cC, \cD)$ is an $A_\infty$ category. For a conceptual explanation, as well as explicit formulas for the differential and composition, cf. \cite{Lu}, \cite{BLM}, \cite{K1}.

Furthermore, for DG categories $\cC$ and $\cD$ there are $A_\infty$ morphisms
\begin{equation}\label{eq:mult fi}
\cC \otimes \Fun_\infty(\cC, \cD)\to \cD
\end{equation}
(the action) and
\begin{equation}\label{eq:mult fi1}
\Fun_\infty(\cD, \cE)\otimes \Fun_\infty(\cC, \cD)\to \Fun_\infty(\cC, \cE)
\end{equation}
(the composition). This follows from the conceptual explanation cited below; in fact these pairing were considered already in \cite{Ko}. As a consequence, there are pairings
\begin{equation}\label{eq:mult fi2}
{\CC}^-_{\bullet}(\cC)\otimes{\CC}^-_{\bullet}( \Fun_\infty(\cC, \cD))\to {\CC}^-_{\bullet}(\cD)
\end{equation}
and
\begin{equation}\label{eq:mult fi3}
{\CC}^-_{\bullet}(\Fun_\infty(\cD, \cE))\otimes{\CC}^-_{\bullet}( \Fun_\infty(\cC, \cD))\to {\CC}^-_{\bullet}(\Fun_\infty(\cC, \cE))
\end{equation}
 Recall that Getzler and Jones constructed an explicit $A_{\infty}$ structure on the negative cyclic complex of an associative commutative algebra. The formulas involve the shuffle product and higher cyclic shuffle products; cf. \cite{GJ}, \cite{L}. When the algebra is not commutative, the same formulas may be written, but they do not satisfy the correct identities. One can, however, define external Getzler-Jones products for algebras and, more generally, for DG categories by the same formulas. One gets maps 
$${\CC}^-_{\bullet}(\cC_1)\otimes \ldots \otimes{\CC}^-_{\bullet}(\cC_n)\to {\CC}^-_{\bullet}(\cC_1\otimes \ldots\cC_n)[2-n]$$
which satisfy the usual $A_\infty$ identities.
To get \eqref{eq:mult fi2} and \eqref{eq:mult fi3}, one combines these products with \eqref{eq:mult fi} and \eqref{eq:mult fi1}.

\begin{example}\label{ex:one functor}
Let $F$ be an $A_\infty$ functor from $\cC$ to $\cD$. Then ${\operatorname{id}}_F$ is a chain of ${\CC}^{-}(\Fun_{\infty}(\cC, \cD))$ (with $n=0$). The pairing \eqref{eq:mult fi2} with this chain amounts to the map of the negative cyclic complexes induced by the $A_\infty$ functor $F$:
$$f_0\otimes \ldots \otimes f_n\mapsto \sum\pm F_{k_0}(\ldots f_0\ldots)\otimes F_{k_1}(\ldots)\otimes \ldots \otimes F_{k_m}(\ldots)$$
The sum is taken over all cyclic permutations of $f_0,\ldots,f_n$ and all $m,\,k_0,\,\ldots, \,k_m$ such that $f_0$ is inside $F_{k_0}$.
\end{example}

\begin{remark}\label{remark:tamark}
The action \eqref{eq:mult fi} and the composition \eqref{eq:mult fi}  are parts of a very nontrivial structure that was studied in \cite{Ta}.
\end{remark}

As a consequence, this gives an $A_\infty$ category structure $\CC^-({\bf \Fun_\infty})$ whose objects are $A_\infty$ categories and whose complexes of morphisms are negative cyclic complexes $\CC^-_{\bullet}(\Fun_\infty(\cD, \cE)).$

>From a less conceptual point of view, pairings \eqref{eq:mult fi2} and \eqref{eq:mult fi3} were defined, in partial cases, in \cite{NT1} and \cite{NT}. The $A_\infty$ structure on $\CC^-({\bf \Fun_\infty})$ was constructed (in the partial case when all $f$ are identity functors) in \cite{TT}. Cf. also \cite{T1} for detailed proofs.
\subsection{The prefibered version}\label{ss:the prefibered version}
We need the following modification of the above constructions. Let $\cB$ be a
category. Consider, instead of a single DG category $\cD$, a family of DG categories $\cD_i$, $i\in\mathrm{Ob}(\cB),$ together with a family of DG functors $f^*:{\cD}_i\leftarrow{\cD}_j,$ $f\in {\cB}(i,j),$ satisfying $(fg)^*=g^*f^*$ for any $f$ and $g$. In this case we define a new DG category $\cD:$
$${\rm{Ob}}(\cD)=\coprod _{i\in {\rm{Ob}}(\cB)}{\rm{Ob}}(\cD_i)$$
and, for $a\in {\rm{Ob}}(\cD_i),$ $b\in {\rm{Ob}}(\cD_j),$
$$\cD(a,b)=\oplus_{f\in \cB(i,j)}\cD_i (a,f^*b).$$
The composition is defined by
$$(\varphi, f)\circ (\psi,g)=(\varphi\circ f^*\psi, f\circ g)$$
for $\varphi \in \cD_i(a,f^*b)$ and $\psi\in \cD_j(b,g^*c).$

We call the DG category $\cD$ {\em a DG category over $\cB,$} or, using the language of \cite{Gil}, {\em a prefibered DG category over $\cB$ with a strict cleavage.} There is a similar construction for $A_\infty$ categories.

Let $\cC,$ $\cD$ be two DG categories over $\cB$. An $A_\infty$ functor $F:\cC\to\cD$ is called {\em an $A_\infty$ functor over $\cB$} if for any $a\in {\rm{Ob}}(\cC_i)$ $Fa\in {\rm{Ob}}(\cD_i),$  and for any $a_k\in {\rm{Ob }}(\cC_{i_k}),$ $(\varphi_k, f_k)\in \cC(a_k, a_{k-1}),$ $k=1,\ldots,n,$
$$F_n((\varphi_n, f_n),\ldots,(\varphi_1, f_1))=(\psi,f_1\ldots f_n)$$
for some $\psi\in \cD_{i_n}.$ One defines a morphism over $\cB$ of two $A_\infty$ functors over $\cB$ by imposing a restriction which is identical to the one above. We get a DG category $\Fun_\infty^{\cB}(\cC, \cD).$ As in the previous section, there are $A_\infty$ functors
\begin{equation}\label{eq:mult fi fi}
\cC \otimes \Fun_\infty^{\cB}(\cC, \cD)\to \cD
\end{equation}
(the action) and
\begin{equation}\label{eq:mult fi fi 1}
\Fun_\infty^{\cB}(\cD, \cE)\otimes \Fun_\infty^{\cB}(\cC, \cD)\to \Fun_\infty^{\cB}(\cC, \cE)
\end{equation}
(the composition), as well as
\begin{equation}\label{eq:mult fi fi 2}
{\CC}^-_{\bullet}(\cC)\otimes{\CC}^-_{\bullet}( \Fun_\infty^{\cB}(\cC, \cD))\to {\CC}^-_{\bullet}(\cD)
\end{equation}
and
\begin{equation}\label{eq:mult fi fi 3}
{\CC}^-_{\bullet}(\Fun_\infty^{\cB}(\cD, \cE))\otimes{\CC}^-_{\bullet}( \Fun_\infty^{\cB}(\cC, \cD))\to {\CC}^-_{\bullet}(\Fun_\infty^{\cB}(\cC, \cE))
\end{equation}
\subsubsection{}\label{sss:proj cat} We need one more generalization of the above constructions. It is not necessary if one adopts the convention from Remark \ref{rmk:simplification}.

Suppose that instead of $\cB$ we have a diagram of categories indexed by a category ${\bf U}$ (in other words, a functor from ${\bf U}$ to the category of categories. In our applications, ${\bf U}$ will be the category of open covers). Instead of a $\cB$-category $\cD$ we will consider a family of $\cB_u$-categories $\cD_u$, $u\in {\operatorname{Ob}}({\bf U}),$ together with a functor $\cD_v\to \cD_u$ for any  morphism $u\to v$ in ${\bf U}$, subject to compatibility conditions that are left to the reader. The inverse limit of categories $\underset{\bf{U}}{\liminv} \cD_u$ is then a category over the inverse limit $\underset{\bf{U}}{\liminv} \cB_u$. We may proceed exactly as above and define the DG category of $A_{\infty}$ functors over ${\liminv} \cB_u$ from ${\liminv} \cD_u$ to ${\liminv} \cE_u$, etc., with the following convention: the space of maps from the inverse product, or from the tensor product of inverse products, is defined to be the inductive limit of spaces of ma!
 ps from (tensor products of) individual constituents.

In this new situation, the pairings \eqref{eq:mult fi fi 1} and \eqref{eq:mult fi fi 3} still exist, while \eqref{eq:mult fi fi 2} turns into

%\begin{equation}\label{eq:mult fi fi 5}
%\Fun_\infty^{\cB}(\cC, \cD)\to\limdir {\underline{Hom}}(\cC_u,\cD)
%\end{equation}
\begin{equation}\label{eq:mult fi fi 6}
{\CC}^-_{\bullet}( \Fun_\infty^{\cB}(\cC, \cD))\to \limdir{\underline{Hom}}({\CC}^-_{\bullet}(\cC_u),\liminv{\CC}^-_{\bullet}(\cD_v))
\end{equation}
\subsection{The trace map for stacks}\label{ss:trace map for stacks}
\subsubsection{From perfect to very strictly perfect complexes}\label{From perfect to very strictly perfect complexes}
Let $M$ be a space with a stack $\cA.$ Consider an open cover ${\mathfrak U}=\{U_i\}_{i\in I}$ such that the stack $\cA$ can be represented by a datum $\cA_i, G_{ij}, c_{ijk}$. Let $\cB_{\mathfrak U}$ be the category whose set of objects is $I$ and where for every two objects $i$ and $j$ there is exactly one morphism $f:i\to j.$ Put $\cC_{\mathfrak U}=k[\cB_{\mathfrak U}],$ i.e. $(\cC_{\mathfrak U})_i=k$ for any object $i$ of $\cB_{\mathfrak U}.$

There is a standard isomorphism of the stack $\cA|{U_i}$ with the trivial stack associated to the sheaf of rings $\cA_i.$ Therefore one can identify twisted modules on $U_i$ with sheaves of $A_i$-modules. We will denote the twisted module corresponding to the free module $\cA_i$ by the same letter $\cA_i.$

\begin{definition}\label{dfn:very strictly perfect} Define the category of {\em very strictly perfect complexes} on any open subset of $U_i$ as follows. Its objects are pairs $(e,d)$ where $e$ is an idempotent endomorphism of degree zero of a free graded module $\sum_{a=1}^N\cA_i[n_a]$ and $d$ is a differential on ${\operatorname{Im}}(e).$ Morphisms between $(e_1,d_1)$ and $(e_2,d_2)$ are the same as morphisms between ${\operatorname{Im}}(e_1)$ and ${\operatorname{Im}}(e_2)$ in the DG category of complexes of modules. A parallel definition can be given for the category of complexes of modules over an associative algebra.
\end{definition}

 Let $(\cD_{\mathfrak U})_i$ be the category of very strictly perfect complexes of twisted $\cA$-modules on $U_i$. By ${\bf U}$ we denote the category of open covers as above.

Strictly speaking, our situation is not exactly a partial case of what was considered in \ref{ss:the prefibered version}. First, $(\cD_{\mathfrak U})_i$ is a presheaf of categories on $U_i$ (in the most naive sense, i.e. it consists of a category $(\cD_{\mathfrak U})_i(U)$ for any $U$ open in $U_i,$ and a functor $G_{UV}:(\cD_{\mathfrak U})_i(V)\to (\cD_{\mathfrak U})_i(U)$ for any $U\subset V,$ such that $G_{UV}G_{VW}=G_{UW}$). Second, $f^*$ are defined as functors on the subset $U_i\cap U_j.$ Also, the pairing \eqref{eq:mult fi fi 2} and its generalization \eqref{eq:mult fi fi 6} are defined in a slightly restricted sense: they put in correspondence to a cyclic chain $i_0\to i_n\to i_{n-1}\to \ldots \to i_0$ a cyclic chain of the category of very strictly perfect complexes of $\cA$-modules on $U_{i_0}\cap\ldots \cap U_{i_n}.$ Finally, in the notation of \ref{sss:proj cat}, for a morphism  $f:{\mathfrak U}\to {\mathfrak V}$ in ${\bf U}$ and an object $j$ of $I_{\mathfrak V}!
 ,$ the functor $({\cD}_{\mathfrak V})_j\to({\cD}_{\mathfrak U})_{f(j)}$ induced by $f$ is defined only on the open subset $V_j.$

We put $\cB=\liminv \cB_{\mathfrak U}$ and $\cD=\liminv \cD_{\mathfrak U}.$

Let $\Perf(\cA)$ be the DG category of perfect complexes of twisted $\cA$-modules on $M$. We denote the sheaf of categories of very strictly perfect complexes on $M$ by $\Perf^{\operatorname{vstr}}(\cA).$ If $Z$ is a closed subset of $M$ then by $\Perf_Z(\cA)$ we denote the DG category of perfect complexes of twisted $\cA$-modules on $M$ which are acyclic outside $Z$.
\begin{definition}\label{dfn:complex of cc of matr} Define
$$
\check{C}^{-\bullet}(M, \CC^-_\bullet(\Matr_\tw(\cA)))=\underset{\mathfrak{U}}{\limdir} \prod_{{\sigma_0}\subset {\sigma_1} \subset \ldots \subset {\sigma_p}}{\CC}^-_{\bullet}(\Matr^{\sigma_p}_\tw(\cA))$$
where $\sigma_i$ run through simplices of ${\mathfrak U}$. The total differential is $b+uB+{\check{\partial}}$ where
$${\check{\partial}}s_{\sigma_0\ldots \sigma_p}=\sum_{k=0}^{p-1}(-1)^k s_{\sigma_0\ldots{ \widehat{\sigma_k}}\ldots \sigma_p}+(-1)^ps_{\sigma_0\ldots \sigma_{p-1}}|U_{\sigma_p} $$
For a closed subset $Z$ of $M$ define
$\check{C}_Z^{-\bullet}(M, {\CC}^-_{\bullet}(\Matr_\tw(\cA)))$ as
$${\operatorname{Cone}}(\check{C}^{-\bullet}(M, {\CC}^-_{\bullet}(\Matr_\tw(\cA)))\to \check{C}^{-\bullet}(M\setminus Z,\, {\CC}^-_{\bullet}(\Matr_\tw(\cA))))[-1] \ .$$
\end{definition}

Let us construct natural morphisms
\begin{equation}\label{eq:Trace map}{\CC}^-_{\bullet}(\Perf(\cA))\to \check{C}^{-\bullet}(M, {\CC}^-_{\bullet}(\Matr_\tw(\cA)))
\end{equation}
\begin{equation}\label{eq:Trace map Z}{\CC}^-_{\bullet}(\Perf_Z(\cA))\to \check{C}_Z^{-\bullet}(M, {\CC}^-_{\bullet}(\Matr_\tw(\cA)))
\end{equation}
First, observe that
the definition of a twisted cochain and Lemma \ref{lemma:perfect complex} can be reformulated as follows.
\begin{lemma}\label{lemma:perfect complex1}
1. A twisting cochain is an $A_\infty$ functor $\cC\to \cD$ over $\cB$ in the sense of \ref{ss:the prefibered version}.

2. There is an $A_\infty$ functor from the DG category of perfect complexes to the DG category $\Fun_{\infty}(\cC,\cD).$
\end{lemma}
The second part of the Lemma together with \eqref{eq:mult fi fi 2} give  morphisms
$$
\CC^-_\bullet(\Perf(\cA))\to \CC^-_\bullet(\Fun^\cB_\infty(\cC,\cD)) \to \underline{\Hom}(\CC^-_\bullet(\cC), \CC^-_\bullet(\cD)) \ .
$$
As mentioned above, the image of this map is the subcomplex of those morphisms that put in correspondence to a cyclic chain $i_0\to i_n\to i_{n-1}\to \ldots \to i_0$ a cyclic chain of the category of very strictly perfect complexes of $\cA$-modules on $U_{i_0}\cap\ldots \cap U_{i_n}$. We therefore get a morphism
$$
{\CC}^-_{\bullet}(\Perf(\cA))\to \check{C}^{-\bullet}(M, {\CC}^-_{\bullet}(\Perf^{\operatorname{vstr}}(\cA)))
$$
Now replace the right hand side by the quasi-isomorphic complex
$$
\underset{\mathfrak{U}}{\limdir}\prod_{{\sigma_0}\subset {\sigma_1} \subset \ldots \subset {\sigma_p}}{\CC}^-_{\bullet}(\Perf^{\operatorname{vstr}}(\cA(U_{\sigma_p})))$$
where $\sigma_i$ run through simplices of ${\mathfrak U}$. There is a natural functor
$$\Perf^{\operatorname{vstr}}(\cA(U_{\sigma_p}))\to \Perf^{\operatorname{vstr}}(\Matr_\tw^{\sigma_p}(\cA))$$
where the right hand side stands for the category of very strictly perfect complexes of modules over the sheaf of rings $\Matr_\tw^{\sigma_p}(\cA)$ on $U_{\sigma_p}.$ This functor acts as follows: to a twisted module $\cM$ it puts in correspondence the direct sum $\oplus_{i\in I_{\sigma_0}}\cM_i$; an element $a_{ij}E_{ij}$ acts via $a_{ij}g_{ij}.$
\subsubsection{From the homology of very strictly perfect complexes to the homology of the algebra}\label{From the homology of very strictly perfect complexes to the homology of the algebra}
Next, let us note that one can replace $ {\CC}^-_{\bullet}(\Perf^{\operatorname{vstr}}(
\Matr_\tw^{\sigma_p}(\cA)))$ by the complex ${\CC}^-_{\bullet}(\Matr_\tw^{\sigma_p}(\cA)):$ indeed, for any associative algebra $A$ there is an explicit trace map
\begin{equation}\label{eq:trace for str per}
 {\CC}^-_{\bullet}(\Perf^{\operatorname{vstr}}(A))\to  {\CC}^-_{\bullet}(A)
\end{equation}
Our construction of the trace map can be regarded as a modification of Keller's argument from \cite{K}. First, recall from \ref{ss:Categories of functors} the internal Getzler-Jones products. The binary product will be denoted by $\times.$ 
 We define the map \eqref{eq:trace for str per} as a composition
$$ {\CC}^-_{\bullet}(\Perf^{\operatorname{vstr}}(A))\to   {\CC}^-_{\bullet}({\operatorname{Proj}}(A))\to   {\CC}^-_{\bullet}({\operatorname{Free}}(A))\to{\CC}^-_{\bullet}(A);$$
the second DG category is the subcategory of complexes with zero differential; the third is the subcategory of complexes of free modules with zero differential.  The morphism on the left is the exponential of the operator $-(1\otimes d)\times?$ opposite to the operator of binary product with the one-chain $1\otimes d$. 
% {\underline{\,\,}}\.\,)$;
The morphism in the middle is ${\operatorname{ch}}(e)\times ?,$ the operator of binary multiplication by the Connes-Karoubi Chern character of an idempotent $e$, cf. \cite{L}. The morphism on the right is the standard trace map from the chain complex of matrices over an algebra to the chain complex of the algebra itself \cite{L}.

Let us explain in which sense do we apply the Getzler-Jones product. To multiply $f_0\otimes \ldots \otimes f_n$ by ${\operatorname{ch}}(e),$ recall that $f_k: \cF_{i_k}\to \cF_{i_{k-1}}$, $\cF_{i_k}$ are free of finite rank, $e_k^2=e_k$ in ${\Hom}(\cF_{i_k},\cF_{i_k}),$ $\cF_{i_{-1}}= \cF_{i_n},$ $e_{-1}= e_n,$ and $f_ke_k=e_{k-1}f_k.$ Write the usual formula for multiplication by ${\operatorname{ch}}(e),$ but, when a factor $e$ stands between $f_i$ and $f_{i+1},$ replace this factor by $e_i.$ Similarly for the morphism on the left: if a factor $d$ stands between $f_i$ and $f_{i+1},$ replace this factor by $d_i$ (the differential on the $i$th module).
This finishes the construction of the morphism \eqref{eq:Trace map}.

Next, we need to refine the map \eqref{eq:trace for str per} as follows.  Recall \cite{D} that for a DG category $\cD$ and for a full DG subcategory $\cD_0$ the DG quotient of $\cD$ by $\cD_0$ is the following DG category. It has same objects as $\cD$; its morphisms are freely generated over $\cD$ by morphisms $\epsilon _i$ of degree $-1$ for any $i\in {\operatorname{Ob}}(\cD_0)$, subject to $d\epsilon _i={\operatorname{id}}_i.$ It is easy to see that the trace map \eqref{eq:trace for str per} extends to the negative cyclic complex of the Drinfeld quotient of
$\Perf^{\operatorname{vstr}}(A)$ by the full DG subcategory of acyclic complexes. Indeed, a morphism in the DG quotient is a linear combination of monomials $f_0\epsilon_{i_0}f_1\epsilon _{i_1}\ldots \epsilon_{i_{n-1}}f_n$ where $f_k: \cF{i_k}\to\cF_{i_{k-1}}$ and $\cF_{i_k}$ are acyclic for $k=0, \ldots, n-1.$ An acyclic very strictly perfect complex is contractible. Choose contracting homotopies $s_k$ for $\cF_{i_k}$. Replace all the monomials $f_0\epsilon_{i_0}f_1\epsilon _{i_1}\ldots \epsilon_{i_{n-1}}f_n$ by $f_0s_0 f_1 s_1\ldots s_{n-1}f_n$. Then apply the above composition to the resulting chain of
 ${\CC}^-_{\bullet}(\Perf^{\operatorname{vstr}}(A))$. We obtain for any associative algebra $A$
 \begin{equation}\label{eq:Trace map loc}{\CC}^-_{\bullet}(\Perf^{\operatorname{vstr}}(A)_{\operatorname{Loc}})\to  {\CC}^-_{\bullet}(A)
\end{equation}
where $_{\operatorname{Loc}}$ stands for the Drinfeld localization with respect to the full subcategory of acyclic complexes.

To construct the Chern character with supports, act as above but define $\cD_i$ to be the Drinfeld quotient of the DG category $\Perf^{\operatorname{vstr}}(\cA(U_i))$ by the full subcategory of acyclic complexes. We get a morphism
$$
{\CC}^-_{\bullet}(\Perf(\cA))\to \check{C}^{-\bullet}(M, {\CC}^-_{\bullet}(\Perf^{\operatorname{vstr}}(\cA)_{\operatorname{Loc}}))\to  \check{C}^{-\bullet}(M, {\CC}^-_{\bullet}(\cA))
$$
From this, and from the fact that the negative cyclic complex of the localization of ${\operatorname{Perf}}_Z$ is canonically contractible outside $Z$, one gets easily the map \eqref{eq:Trace map Z}.
 %${\CC}^_{\bullet}(\Perf^{\operatorname{vstr}}(A))$.
 \subsection{Chern character for stacks}\label{ss:Chern character for stacks}
%Now observe that there is a natural Chern character
%\begin{equation}
%\label{eq:McCarthy}
%K_\bullet(\Perf(\cA))\to {\operatorname{HC}}^-_{\bullet}(\Perf(\cA)_{\operatorname{Loc}})
%\end{equation}
%whose composition with \eqref{eq:Trace map loc}
Now let us construct the Chern character
\begin{equation}
\label{eq:ch}
K_\bullet(\Perf(\cA))\to \check{{\mathbb H}}^{-\bullet}(M, {\CC}^-_{\bullet}(\Matr_\tw(\cA)))
\end{equation}
\begin{equation}
\label{eq:ch Z}
K_\bullet(\Perf_Z(\cA))\to \check{{\mathbb H}}_Z^{-\bullet}(M, {\CC}^-_{\bullet}(\Matr_\tw(\cA)))
\end{equation}
First, note that the $K$ theory in the left hand side can be defined as in \cite{TV}; one can easily deduce from \cite{MC} and \cite{K2}, section 1, the Chern character from $K_\bullet(\Perf(\cA))$ to the homology of the complex ${\operatorname{Cone}}({\CC}^-_{\bullet}(\Perf_{\operatorname{ac}}(\cA))\to {\CC}^-_{\bullet}(\Perf(\cA)))$. Here $\Perf_{\operatorname{ac}}$ stands for the category of acyclic perfect complexes.

Compose this Chern character  with the trace map of \ref{ss:trace map for stacks}. We get a Chern character from $K_\bullet(\Perf(\cA))$ to
$${\check{\mathbb H}}^{-\bullet} (M,{\operatorname{Cone}}({\CC}^-_{\bullet}(\Perf^{\operatorname{vstr}}_{\operatorname{ac}}(\cA)_{\operatorname{Loc}})\to {\CC}^-_{\bullet}(\Perf^{\operatorname{vstr}}(\cA)_{\operatorname{Loc}})))$$
One gets the Chern characters \eqref{eq:ch}, \eqref{eq:ch Z} easily by combining the above with \eqref{eq:Trace map loc}.
\subsection{The case of a gerbe}
If $\cA$ is a gerbe on $M$ corresponding to a class $c$ in $H^2(M, \cO _M^*),$ then (in the $C^{\infty}$ case) the right hand side of \eqref{eq:ch} is the cohomology of $M$ with coefficients in the complex of sheaves $$\Omega ^{-\bullet}[[u]], ud_{\operatorname {DR}}+u^2 H\wedge$$ where $H$ is a closed three-form representing the three-class of the gerbe. In the holomorphic case, the right hand side of \eqref{eq:ch} is computed by the complex $\Omega ^{-\bullet, \bullet}[[u]], {\overline{\partial}}+\alpha \wedge + u\partial$ where $\alpha$ is a closed $(2,1)$ form representing the cohomology class $\partial{\operatorname{log}}c.$ This can be shown along the lines of \cite{BGNT}, Theorem 7.1.2.

\end{document}